\documentstyle[12pt]{article}

\begin{document}
\newcommand{\nt}{\noindent}

\begin{center}

{\bf Quantum Schubert Calculus}

\medskip

Aaron Bertram\footnote{Partially supported by NSF grant DMS-9218215
and a Sloan fellowship}

\end{center}

\nt {\bf 1. Introduction:} 
Classical Schubert calculus, by which I mean the 
formulas of Giambelli and Pieri encoding the product structure of
the cohomology ring of a complex Grassmannian, has been an essential 
tool in enumerative algebraic geometry for over a century.  

\medskip 

String theorists (notably Witten \cite{W}) recently introduced the notion
of a ``quantum'' deformation of the cohomology ring of a smooth
projective variety $X$. This quantum deformation, or quantum cohomology ring, 
as it is often called, is an algebra over a formal-power-series
ring which specializes to the ordinary cohomology ring, and which is 
defined in terms of intersection data (the Gromov-Witten
invariants) on all the spaces
of holomorphic maps from pointed curves of genus zero to $X$.

\medskip
 
A rigorous definition of the Gromov-Witten invariants, together with a
verification of the algebra structure of 
these quantum deformations, has been established by 
two schools, namely the symplectic school of Ruan-Tian
(\cite {RT}) and the algebro-geometric school of 
Kontsevich-Manin (\cite{KM}). One interesting variant of 
the quantum deformation is a ``small'' deformation of the 
cohomology ring (terminology taken from \cite{F2}) which 
is an algebra over a polynomial ring (hence of finite-type
over {\bf C}) sitting between the full quantum deformation
and the cohomology ring itself. This ``small'' quantum cohomology ring 
can be defined independently, and we will do so in the 
Grassmannian case, where it turns out to be an algebra 
over a polynomial ring in one variable. We'll
let $q$ stand for the variable.

\medskip
 
In this paper, the rules for the Schubert calculus are modified so 
that they are valid in the small quantum cohomology ring.
In other words, whereas the Giambelli and Pieri formulas
are valid in the cohomology ring of a Grassmannian, higher order terms
(in $q$) may appear when the corresponding products are 
taken in this ring.
The main
result here is the computation of these higher order terms.   
Our computation relies on the recursive properties of a particular
smooth compactification (the Grothendieck quot scheme)
of the space of holomorphic maps of a fixed 
degree from ${\bf P}^1$ to a Grassmannian.
 
\bigskip

In order to fix notation and refresh the reader's memory, 
we begin with an overview (following \cite{GH})
of the classical Schubert calculus before continuing with the introduction.

\begin{tabbing}

Let: \=  $V$ be a vector space over ${\bf C}$ of dimension $n$, \\ 
\> $0 = V_0 \subset V_1 \subset ... \subset V_n = V$ be a full flag for $V$,\\ 
\> $G := G(n-k,n)$ be the Grassmannian of $n-k$-dim'l subspaces of $V$, \\  
\> $\Lambda_x \subset V$ be the subspace corresponding to a point $x\in G$.  

\end{tabbing}

Given an $n-k$-tuple of integers $\vec a := (a_1,...,a_{n-k})$ satisfying  
the inequalities $k \ge a_1 \ge ... \ge a_{n-k} \ge 0$, 
let: 
$$W_{\vec a} = \{ x \in G | \ \mbox{dim}(\Lambda_x\cap V_{k+i-a_i}) \ge i\}.$$

Then $W_{\vec a}$ is a subvariety of $G$ of complex codimension $|\vec a|
:= \sum _{i=1}^{n-k} a_i$. Let
$$\sigma_{\vec a} \in \ \mbox H ^{2|\vec a|}(G,{\bf C})$$ 
be the corresponding element in cohomology.

\medskip

One calls $W_{\vec a}$ 
the Schubert variety associated to $\vec a$ (and the given flag). 
The cohomology classes $\sigma_{\vec a}$
produce a vector-space basis for $\mbox H^*(G,{\bf C})$ as the $\vec a = (a_1,...,a_{n-k})$
range over all tuples of integers with the given constraints. 

\medskip

The  Schubert varieties $W_a := W_{(a,0,...,0)}$ are called 
special Schubert varieties. The corresponding cohomology classes 
$\sigma_a$ coincide with the image in cohomology of the chern classes 
$c_a(Q)$ where $Q$ is 
the universal quotient bundle on $G$. 
These special cohomology classes generate the cohomology 
ring of the Grassmannian as an {\bf algebra} over ${\bf C}$ via
the following determinantal formula:

\medskip

\nt {\bf Giambelli's Formula:} By convention, let $\sigma_a = 0$ if $a < 0$
 or $a > k$. Then:
$$\sigma_{\vec a} = \Delta_{\vec a}(\sigma _*) := \left|\begin{array}{ccccc}
\sigma_{a_1} & \sigma_{a_1+1} & \sigma_{a_1+2} & \cdots & \sigma_{a_1+n-k-1} \\
\sigma_{a_2-1} & \sigma_{a_2} & \sigma_{a_2+1} & \cdots & \sigma_{a_2+n-k-2} \\
\sigma_{a_3-2} & \sigma_{a_3-1} & \sigma_{a_3}  \\
\vdots & & & & \vdots \\
\sigma_{a_{n-k} - (n-k) + 1} & & \cdots & & \sigma_{a_{n-k}}\end{array}\right|$$

\bigskip

The other ``Italian'' formula explicitly computes the intersections  of a 
special cohomolgy class and a general one:

\bigskip

\nt {\bf Pieri's Formula:} 
The product of $\sigma_a$ and $\sigma _{\vec a}$ in $\mbox H^*(G,{\bf C})$
is:
$$\sigma_a \cdot \sigma_{\vec a} = p_{a,\vec a}(\sigma_{\vec *}) := 
\sum_{\vec b}\sigma_{\vec b}$$
where the $\vec b$ vary over all $n-k$-tuples satisfying:

\medskip

$|\vec b| = a + |\vec a|$ and $k \ge b_1 \ge a_1 \ge ... \ge b_{n-k} \ge a_{n-k} \ge 0$.

\bigskip

The two Italian formulas determine the ring 
structure on the cohomology ring $\mbox H^*(G,{\bf C})$. Indeed, (see \cite{BT},
Proposition 23.2)
one obtains a convenient presentation of the cohomology ring as:
$$(*):{\bf C}[X_1,...,X_k]/(Y_{n-k+1}(X_*),...,Y_n(X_*)) \cong \ \mbox H^*(G,{\bf C}), \ \ \ 
X_a \mapsto \sigma _a$$
where $Y_i(X_*)$ is the coefficient of $t^i$ in the formal-power-series inverse
of the polynomial $1+X_1t+ ...+X_kt^k$.

\bigskip

Finally, recall that if we set $\vec a^c := (k-a_{n-k},k-a_{n-k-1},...,k-a_1)$,
then $\sigma_{\vec a}$ and $\sigma_{\vec a^c}$ are 
Poincar\'e dual. Equivalently, if we let 
$\langle W_{\vec a},W_{\vec b}\rangle$ denote the intersection 
number in $G$ of general translates of $W_{\vec a}$ and $W_{\vec b}$ 
(which is set to zero
if $|\vec a| + |\vec b| \ne \ \mbox{dim}(G) = k(n-k)$), then 
$$\langle W_{\vec a},W_{\vec b}\rangle = \left\{
\begin{array}{ll} 1 \ \ & \mbox {if} \ \vec b = \vec a^c \\
0  & \mbox{otherwise}\end{array}\right.$$

\bigskip

The multiplication on $\mbox H^*(G,{\bf C})$ can be understood solely 
in terms of intersection numbers as follows. If $\sigma_{\vec a_1},...,
\sigma_{\vec a_N}$ are cohomology classes corresponding to Schubert 
varieties, let $\langle W_{\vec a_1},...,W_{\vec a_N}\rangle$ be
the intersection number of general translates of the Schubert varieties,
set to zero, as above, if $\sum_{i=1}^N|\vec a_i| \ne \ \mbox{dim}(G)$. Then
in the cohomology ring of $G$,
$$\sigma_{\vec a_1}\cdot ...\cdot \sigma_{\vec a_N} = 
\sum_{\vec a} \langle W_{\vec a}, W_{\vec a_1},...,W_{\vec a_N}\rangle \ \sigma_{\vec a^c}.$$

This is an immediate consequence of the fact that the cohomology
classes $\sigma_{\vec a}$ satisfy the Poincar\'e duality 
property above.

\bigskip

The key idea behind the quantum deformations is to introduce 
``higher order terms'' into the product by considering a sequence of
intersection numbers, starting with the intersections on $G$ itself. 
To be more precise, here is a definition
for the small quantum deformation. 

\bigskip

For each integer $d \ge 0$, let 
$\langle W_{\vec a_1},... W_{\vec a_N}\rangle_d \in {\bf Z}$
be the ``Gromov-Witten'' intersection number defined as follows. Choose general
points $p_1,...,p_N \in {\bf P}^1$ and general translates of 
the $W_{\vec a_i}$. Then 
$\langle W_{\vec a_1},...,W_{\vec a_N}\rangle_d$ is ``by definition'' the
number of holomorphic maps $f:{\bf P}^1 \rightarrow G$ of degree $d$  
with the 
property that $f(p_i) \in W_{\vec a_i}$ for all $i = 1,...,N$
(and zero if the sum of the $|\vec a_i|$ is such that one expects this 
number not to be finite).  In \S 2, we make rigorous sense out of
this definition and reinterpret it as an intersection of generalized
Schubert cohomology classes in a Grothendieck quot scheme, which 
in this case happens to be a smooth, projective variety of dimension
$nd + \ \mbox{dim}(G)$. 
Notice that in particular, 
the Gromov-Witten number $\langle W_{\vec a_1},... W_{\vec a_N}\rangle_0$
is the original intersection number in $G$.

\medskip

The ``quantum'' product of $\sigma_{\vec a_1},...,\sigma_{\vec a_N}$,
which we will denote with asterisks as 
$\sigma_{\vec a_1}*...*\sigma_{\vec a_N}$,
is defined as follows. Let $q$ be a formal variable. Then: 
$$\sigma_{\vec a_1}*...*\sigma_{\vec a_N} = 
\sum_{d \ge 0} q^d (\sum_{\vec a}  \langle W_{\vec a}, W_{{\vec a}_1},...,W_{{\vec a}_N}\rangle _d
\ \sigma_{\vec a^c})$$

Notice that setting $q= 0$, one recovers the original product. Notice also
that this sum is finite, because the dimensions of the spaces of 
holomorphic maps from ${\bf P}^1$ to $G$ 
increase with $d$.

\bigskip

The really surprising aspect of quantum cohomology is the following: 

\medskip

\nt {\bf Associativity Theorem:} Extend the quantum product to 
a product on elements of $H^*(G,{\bf C})[q]$ by linearity and by setting
$$(\sigma_{\vec a_1}q^{d_1})* ... * (\sigma_{\vec a_N}q^{d_N})
= (\sigma_{\vec a_1}*...*\sigma_{\vec a_N})q^{(d_1 + ... + d_N)}.$$
Then the pairwise quantum product is associative and 
gives H$^*(G,{\bf C})[q]$ the structure
of a ${\bf C}[q]$-algebra. The quantum product of more than
two terms agrees with the product
in this ring.

\bigskip

It would be confusing to refer to this ring as $H^*(G,{\bf C})[q]$
because the quantum product is not the same as the natural product on this 
polynomial ring, so we make the following:

\medskip

\nt {\bf Definition (of the small quantum ring):} 
The small quantum cohomology ring $QH^*(G)$ is by definition 
the vector space H$^*(G,{\bf C})[q]$ equipped with the 
extended quantum product. 

\bigskip

As I said earlier, this theorem is a special case of more general 
associativity results in 
Ruan-Tian \cite{RT} or Kontsevich-Manin \cite{KM}. It is also a very powerful
theorem, as it tells us that all quantum products are 
determined by the pairwise products(!) For example, Siebert and Tian (\cite{ST}),
following ideas of Witten, use this idea to reduce the proof of the following 
presentation for $QH^*(G)$ to a single computation for degree one maps:
$$(*)_q:{\bf C}[X_1,...,X_k,q]/(Y_{n-k+1}(X_*),...,Y_n(X_*) - (-1)^{k-1}q) \cong \ QH^*(G), \ \ \ 
$$
where $X_a \mapsto \sigma_a, q\mapsto q$ and the $Y_i(X_*)$ are defined as in 
$(*)$.

\medskip

In this paper, we will compute
versions of the Italian formulas where the ordinary multiplication  
is replaced by quantum multiplication. Unlike the presentation for 
the quantum cohomology ring above, it seems that the best way to 
approach this problem is not by invoking the
associativity theorem, even though the quantum product is, of course,
determined by $(*)_q$. Rather, 
both formulas follow rather quickly from a theorem of Kempf and Laksov
(\cite{KL}) once we have analyzed the relevant Grothendieck quot scheme in \S 3. 
To be precise, we will prove the following formulas in \S 4:

\medskip

\nt {\bf Quantum Giambelli:} 
$$\sigma_{\vec a} = \Delta_{\vec a}(\sigma_*),$$ 
when the determinant is evaluated in $QH^*(G)$ using the 
quantum product.
In other words, no higher order terms arise from the Giambelli determinant! 

\bigskip

\nt {\bf Quantum Pieri:} 
$$\sigma_a * \sigma_{\vec a} = 
p_{a,\vec a}(\sigma_{\vec *}) + 
q(\sum_{\vec c}\sigma_{\vec c})$$
where the $\vec c$ range over all $n-k$-tuples
satisfying:
$$|\vec c| = a + |\vec a| - n\ \mbox{and} \  
a_1-1 \ge c_1 \ge a_2-1\ge ... \ge a_{n-k}-1
\ge c_{n-k} \ge 0.$$

Notice that as is the case with the classical Giambelli and Pieri formulas, 
the quantum versions determine all the quantum products.

\medskip

In \S 5, as a quick application of quantum Giambelli, we see
that a residue formula of Vafa and Intriligator computing the 
Gromov-Witten numbers for special Schubert varieties can readily 
be modified to compute all the Gromov-Witten numbers.

\bigskip

\nt {\bf Final Remarks:} By substituting the Giambelli determinant,
one of course has the identity:
$\sigma_{\vec a} * \sigma_{\vec a_1}*...*\sigma_{\vec a_N} = 
\Delta_{\vec a}(\sigma_*) * \sigma_{\vec a_1}*...*\sigma_{\vec a_N}$
in $QH^*(G)$ for any Schubert cohomology classes
$\sigma_{\vec a}$ and $\sigma_{\vec a_1},...,\sigma_{\vec a_N}$. 
Similarly one can substitute for a product $\sigma_a*\sigma_{\vec a}$
using quantum Pieri. This is obvious once the associativity theorem
is established. But it can be (and was orignally) proved
directly using the methods of this paper without appealing to 
the quantum ring, and can indeed be used to obtain an 
independent proof of the associativity theorem in this context. 
(Quantum Giambelli
implies H$^*(G,{\bf C})[q]$ is a quotient of the polynomial
ring ${\bf C}[x_1,...,x_k,q]$ and quantum Pieri implies that the 
kernel is an ideal, putting the quotient ring structure, which is 
the quantum product, on H$^*(G,{\bf C})[q]$.)

\medskip

Finally, it is possible to derive quantum Pieri from  
quantum Giambelli and the presentation $(*)_q$
of the quantum cohomology ring, bypassing the delicate geometric 
arguments in \S 4. Quantum Giambelli itself, however, seems
to require a special proof.

\bigskip

\nt {\bf Acknowledgements:} 
I was introduced to quantum cohomology and the formula of Vafa 
and Intriligator by Richard 
Wentworth (see \cite{BDW} or \cite{R} for a completely different way of 
computing the intersection numbers!). I also profited from 
conversations with Youngbin Ruan, Henri Shahrouz and 
Ionut Ciocan-Fontanine on the papers 
mentioned above, as well as 
many much-appreciated comments from Bill Fulton.

Finally, I would like to thank Lowell Abrams, who pointed out 
an error in an earlier formulation of quantum Pieri. While 
correcting that error, I decided to revise
the paper and ``quantize'' the title, which was formerly ``Modular
Schubert Calculus.''

\newpage

\nt {\bf 2. Intersections on the Space of Maps:} In this section, we  will
rigorously define the intersection of Schubert varieties on the moduli spaces 
${\cal M}_d$ of holomorphic maps of degree $d$ from 
${\bf P}^1$ to $G$.
We will first
prove a moving lemma, stating that the Schubert varieties can be made to intersect
in points when they ought to. Then we will prove a cohomological lemma, 
interpreting the intersection
as the intersection of cohomology classes in a given
smooth, {\it projective} variety.  

\bigskip

We begin with ${\cal M}_d$ itself. The usual way to prove
that the space of maps is represented by a quasiprojective scheme
is to embed it as an open set in a Hilbert scheme based on the product
${\bf P}^1\times G$. However, there is another moduli space 
available which contains ${\cal M}_d$ as an open 
subscheme, namely Grothendieck's quot scheme, which
will be our compactification of choice. 

\medskip

Recall that a map
$f:{\bf P}^1 \rightarrow G$ is equivalent to specifying 
a quotient vector bundle $V\otimes {\cal O}_{{\bf P}^1} \rightarrow Q$,
or dually, a subbundle 
$S = Q^* \hookrightarrow V^*\otimes {\cal O}_{{\bf P}^1}$
(modulo
the action of $\mbox {GL}(V)$) where $S$ has degree $-d$ and rank $k$.
The quot scheme will parametrize maps 
$S\hookrightarrow V^*\otimes {\cal O}_{{\bf P}^1}$ that are injective as 
maps of {\it sheaves}. In other words, the cokernel ${\cal F}$ of such a map 
is not
required to be a vector bundle. Specifically, let $\chi(m) = (m+1)(n-k) + d$. 
That is, $\chi$ is the
Hilbert polynomial of a vector bundle of rank $n-k$ and degree $d$ on
${\bf P}^1$. Then:

\medskip 

\nt {\bf Grothendieck's Theorem:} The functor $Quot_\chi(V^*/{\bf P}^1)$
parametrizing flat families of quotients 
$V^*\otimes {\cal O}_{{\bf P}^1} \rightarrow {\cal F}$ of Hilbert polynomial $\chi$
is  representable by a smooth projective variety of dimension
$nd + \mbox{dim}(G)$. 
Moreover, if we let ${\cal Q}_d$
denote this fine moduli space, then 
${\cal M}_d$ is an open subscheme of ${\cal Q}_d$
via the canonical inclusion.

\bigskip

We will call ${\cal Q}_d$ the quot scheme 
compactification of ${\cal M}_d$. (\cite{Gr}
is the standard reference for the proof of Grothendieck's Theorem.) 

\bigskip

Since ${\cal Q}_d$ is a fine moduli space, there is
by definition a universal exact sequence:
$$0 \rightarrow {\cal S}_d \rightarrow V^*\otimes {\cal O}
\rightarrow {\cal T}_d \rightarrow 0,$$
on ${\bf P}^1\times {\cal Q}_d$
which is flat over ${\cal Q}_d$.
The sheaf ${\cal T}_d$ is certainly not usually a vector bundle,
but it is an easy consequence of flatness that the {\it kernel}
${\cal S}_d$ is a vector bundle.

\medskip

We will be very interested in the dual (nonsurjective!) universal map:
$$u: V\otimes {\cal O} \rightarrow {\cal S}_d^*$$

Next, we define the pull-back of Schubert varieties to ${\cal M}_d$.

\medskip

\nt {\bf Definition 2.1:} If $p\in {\bf P}^1$ and $W_{\vec a} \subset G$ is a Schubert variety,
then let:
$$W_{\vec a}(p) = \{ f\in {\cal M}_d\ | \ f(p) \in W_{\vec a} \}$$

We  put a scheme structure on $W_{\vec a}(p)$ via the 
universal evaluation map $ev:{\bf P}^1 \times {\cal M}_d
\rightarrow G$ by redefining:
$$W_{\vec a}(p) := ev^{-1}(W_{\vec a}) \cap 
\left\{ \{p\}\times {\cal M}_d\right\}.$$

We may extend $W_{\vec a}(p)$ as a degeneracy locus 
to the quot scheme:

\medskip

\nt {\bf Definition 2.1A:} Recall the flag of subspaces $0 = V_0 \subset V_1 
\subset ... \subset V_n = V$. 
For each $i=1,...,n-k$,
let $D_{i,a_i} \subset {\bf P}^1 \times {\cal Q}_d$ be the 
largest subscheme 
on which the dimension of the kernel of $u:V_{k-i+a_i} \otimes {\cal O} \rightarrow
{\cal S}^*_d$ is at least $i$, and let $D_{i,a_i}(p)$ be the 
intersection: 
$D_{i,a_i} \cap \left\{ \{p\}\times {\cal Q}_d\right\}$
thought of as a subscheme of ${\cal Q}_d$.
Then exactly as in the definition of $W_{\vec a}$, we define:
$$\overline W_{\vec a}(p) := D_{1,a_1}(p) \cap ... \cap D_{n-k,a_{n-k}}(p).$$

Suppose now that $\vec a_1,...,\vec a_N$ are $(n-k)$-tuples describing
Schubert varieties, and let $A = \sum_{j=1}^{N} |\vec {a_j}|$. Then:

\bigskip

\nt {\bf Moving Lemma 2.2:} For any points $p_1,...,p_N \in {\bf P}^1$,
corresponding general translates of the $W_{{\vec a}_j} \subset G$, 
and a fixed subvariety $Z\subset {\cal M}_d$, the intersection:
$W_{{\vec a}_1}(p_1) \cap ... \cap W_{{\vec a}_N}(p_N) \cap Z$
is either empty, or has pure codimension $A$ in $Z$.

\medskip

{\bf Proof:} It suffices by induction to prove that given 
a subvariety $Z \subset {\cal M}_d$, a point $p\in {\bf P}^1$,
and a general
translate of $W_{\vec a}$, the intersection $Z\cap W_{\vec a}(p)$ is 
empty or has codimension $|\vec a|$ in $Z$. But if we let 
$T \subset G$ be the image of $Z \subset \left\{ \{p\} \times 
{\cal M}_d \right\}$ under the evaluation map $ev$, then 
by an argument of Kleiman (see \cite{H}, III.10.8),
a general translate of $W_{\vec a}$ intersects $T$ in codimension
$|\vec a|$. More generally,  the general translate intersects each locus in $T$ over
which $ev|_Z$ has constant fiber dimension in codimension $|\vec a|$.
Since $W_{\vec a}(p) \cap Z = ev^{-1}(W_{\vec a} \cap T) \cap Z$,
the lemma follows.

\medskip

The lemma implies that when $A = \ \mbox{dim}({\cal M}_d)$,
the schemes $W_{{\vec a}_i}(p_i)$ can be chosen to 
intersect in (reduced) points. Recall that the Gromov-Witten intersection
number is defined as the number of these points when the $p_i$ are in 
general position. However,
it turns out that it suffices for them to be distinct:

\bigskip

\nt {\bf Moving Lemma 2.2A:} If $p_1,...,p_N \in {\bf P}^1$ are
{\it distinct} points, 
then for general choices of $N$ full flags on $V$, an intersection
$\overline W_{{\vec a}_1}(p_1) \cap ... \cap \overline W_{{\vec a}_N}(p_N)$
of generalized Schubert varieties is either empty, 
or has pure codimension $A$ in ${\cal Q}_d$. 
Moreover, the intersection
$W_{{\vec a}_1}(p_1) \cap ... \cap W_{{\vec a}_N}(p_N)$ is Zariski dense in
$\overline W_{{\vec a}_1}(p_1) \cap ... \cap \overline W_{{\vec a}_N}(p_N)$.
In particular,
if $A = \ \mbox{dim}({\cal Q}_d)$, then
$$\overline W_{{\vec a}_1}(p_1) \cap ... \cap \overline W_{{\vec a}_N}(p_N)
= W_{{\vec a}_1}(p_1) \cap ... \cap W_{{\vec a}_N}(p_N).$$

In order to prove this lemma, we will need to analyze the 
structure of the boundary 
${\cal B}_d := {\cal Q}_d - {\cal M}_d$.
This we will do in the next section. For now, we list the main
consequences of the lemma.

\bigskip

\nt {\bf Corollary 2.3:} The cohomology class  
$\sigma_{\vec a} \in \mbox{H}^{2|\vec a|}({\cal Q}_d,{\bf C})$ 
associated to $\overline W_{\vec a}(p)$ is independent of
$p\in {\bf P}^1$ and the choice of flag on $V$. We will call this the
generalized cohomology class associated to the Schubert 
variety $W_{\vec a}$.

\medskip

{\bf Proof:} 
It follows immediately from Definition 2.1A that the 
$\overline W_{\vec a}(p)$ are fibers of a morphism from a subscheme
$X\subset {\bf P}^1\times F \times {\cal Q}_d$ to
${\bf P}^1\times F$, where $F$ is the full flag variety associated
to $V$. 

Because the automorphism groups of ${\bf P}^1$ and $F$ are transitive,
the map $X\rightarrow {\bf P}^1\times F$ is even a fiber bundle, 
of fiber codimension $A$ by the Lemma, and the Corollary follows.
(The reader may check that $X$ is not empty!)

\bigskip

\nt {\bf Corollary 2.4:} If $A = \mbox{dim}({\cal M}_d)$,
then the total degree of the intersection in Lemma 2.2 is independent
of the (general) translates of the $W_{{\vec a}_j}$ and the points $p_j\in {\bf P}^1$
as long as the $p_j$ are distinct.

\medskip

{\bf Proof:} If the $p_j$ are distinct, then the latter part of Lemma
2.2A applies and the intersection number may be interpreted 
as the degree of the product of the $\sigma_{{\vec a}_j}$
in the cohomology ring of the quot scheme.

\medskip

It is very important that the points are distinct in Corollary 2.4.
If the Corollary were true for any collection of points, then the 
quantum Schubert calculus would be trivial! (But see quantum Giambelli in \S 4.)

\bigskip

\nt {\bf Conclusion:} If $\sum_{i=1}^N |\vec a_i| = \ 
\mbox{dim}({\cal M}_d)$,
then the Gromov-Witten number 
$\langle W_{{\vec a}_1},...,W_{{\vec a}_N} \rangle_d$ from the introduction
is well-defined, and coincides with the degree of the product of 
the generalized Schubert cohomology classes $\sigma_{\vec a_1},...,\sigma_{\vec a_N}$ 
in the cohomology ring
of ${\cal Q}_d$. 
 
\medskip

\nt {\bf Remark:} In \S 4, we will use this conclusion to extend the 
definition of the Gromov-Witten numbers to situations where the 
cohomological interpretation is the correct one, and the naive 
definition from \S 1 is not correct. 

\medskip

Thus, the definition of the quantum product is now secure, and the Gromov-Witten
numbers have a cohomological interpretation. 
As another application of the moving lemmas, we prove 
the following formula  for the 
``trivial'' quantum product.

\medskip

\nt {\bf Lemma 2.5:} 
$\sigma_{\vec b} = \sum_{d \ge 0} q^d (\sum_{\vec a}  \langle W_{\vec a}, 
W_{\vec b}\rangle _d
\ \sigma_{\vec a^c})$

\medskip

{\bf Proof:} 
We need to prove that $ \langle W_{\vec a}, W_{\vec b} \rangle _d = 0$
for all pairs of Schubert varieties $W_{\vec a}$ and
$W_{\vec b}$, and all {\it positive} values of $d$.

\medskip 

The Gromov-Witten number is zero if
$|{\vec a}| + |{\vec b}| \ne \ 
\mbox{dim}({\cal M}_d)$ by definition. So we assume equality, and
by the Moving Lemmas, we know that the intersection number is 
realized as the degree of 
$\overline W_{\vec a}(p) \cap \overline W_{\vec b}(o)$ for 
distinct points $o,p \in {\bf P}^1$, and general translates of 
$W_{\vec a}$ and $W_{\vec b}$. Moreover, we know that the intersection
is contained in ${\cal M}_d$.

Now suppose that the intersection is nonempty. Then we have just seen that
there is a 
map $f: {\bf P}^1 \rightarrow G$ of degree $d$ such that
$f(p) \in W_{\vec a}$ and $f(o)\in W_{\vec b}$. But there are 
an entire ${\bf C}^*$ of automorphisms 
$\lambda :{\bf P}^1 \rightarrow {\bf P}^1$ which fix $p,o$, and
the compositions $f\circ \lambda$ all produce {\it different}
elements of $\overline W_{\vec a}(p) \cap \overline W_{\vec b}(o)$.
Since the intersection was proven to be finite, we get a contradiction.

\medskip

(Notice that there is no contradiction in case $d=0$ because if $f$
is a {\it constant} map, then all the $f \circ \lambda$ are the 
same!)

\newpage 

\nt {\bf 3. The Recursive Structure of the Quot Scheme:}
Recall from \S 2 the definition of the boundary of ${\cal Q}_d$:
$${\cal B}_d := {\cal Q}_d
- {\cal M}_d.$$

As we noted in \S 2, the universal quotient sheaf ${\cal T}_d$ 
on ${\bf P}^1\times {\cal Q}_d$ is not locally free. 
In fact, ${\cal M}_d$ is the largest subset $U$ of the quot
scheme with the property that ${\cal T}_d$ has constant rank $n-k$
on ${\bf P}^1\times U$. In the following
theorem, we obtain precise information about a 
stratification of the boundary determined by the loci where 
${\cal T}_d$ has rank at least $n-k+r$.
 
\bigskip

\nt {\bf Theorem 3.1 (Structure Theorem for the Quot Scheme):} 

For all positive integers $r \le k$, let
$\pi_r: {\cal G}_{d,r} \rightarrow 
{\bf P}^1\times {\cal Q}_{d-r}$ be the Grassmann bundle of 
$r$-dimensional quotients of ${\cal S}_{d-r}$
on ${\bf P}^1\times {\cal Q}_{d-r}$, and let 
$\hat {\cal S}_{d-r}$ be the kernel of the tautological quotient 
$\pi_r^*{\cal S}_{d-r} \rightarrow Q$. Then
there are maps $\beta _r:{\cal G}_{d,r} \rightarrow
{\cal Q}_d$ satisfying:

\medskip

(i) If ${\cal T}_d$ has rank at least $n-k+r$ at a point $(p,x)
\in {\bf P}^1\times {\cal Q}_d$, then $x$ is in the image of 
$\beta _r$.

\medskip

(ii) The restriction of $\beta _r$ to 
$\pi_r^{-1}({\bf P}^1\times {\cal M}_{d-r})$
is an embedding.

\medskip

(iii) The preimage of Schubert varieties in ${\cal G}_{d,r}$ is given by
$$\beta _r^{-1}(\overline W_{\vec a}(p))  =
\pi_r^{-1}({\bf P}^1\times \overline W_{\vec a}(p)) \cup 
\widehat W_{\vec a - \vec r}(p)$$
where $\vec r = (r,r,...,r)$, $\widehat W_{\vec a - \vec r}(p) = \cap _{i=1}^{n-k} \hat D_{i,a_i-r}(p)$
and $\hat D_{i,a}(p)$ is the degeneracy locus inside
$\pi_r^{-1}(p\times {\cal Q}_{d-r})$
where the kernel of 
$V_{k-r+i-a}\otimes {\cal O} \rightarrow \hat {\cal S}_{d-r}^*$ has rank
at least $i$.

\bigskip

{\bf Proof of the Structure Theorem:} To construct the maps $\beta _r$, we
need to find bundles ${\cal E}_{d,r} \hookrightarrow V^*\otimes {\cal O}$
on ${\bf P}^1 \times {\cal G}_{d,r}$ which have rank $k$ and relative
degree $-d$ over ${\cal G}_{d,r}$. We obtain these from the
$\pi_r^*{\cal S}_{d-r}$ by 
elementary modifications. Namely, let $\pi_\Delta^* {\cal S}_{d-r} 
\rightarrow \pi_\Delta^*Q$ be the 
pull-back of the tautological
quotient to the preimage of $\Delta
\times {\cal Q}_{d-r}$ in ${\bf P}^1\times {\cal G}_{d,r}$. ($\Delta \subset 
{\bf P}^1\times {\bf P}^1$ is the diagonal.) Let $\pi^*{\cal S}_{d-r}$ be the 
pull-back
of ${\cal S}_{d-r}$ to ${\bf P}^1\times {\cal G}_{d,r}$, and consider the composition:
$$f_{d,r}: \pi^*{\cal S}_{d-r} \rightarrow \pi_{\Delta}^*{\cal S}_{d-r} \rightarrow 
\pi_{\Delta}^*Q.$$
Since the quotient is a vector bundle of rank $r$ supported on a divisor
which intersects each fiber of the projection 
${\bf P}^1 \times {\cal G}_{d,r} \rightarrow {\cal G}_{d,r}$ in a point, 
the kernel of $f_{d,r}$ is a vector bundle ${\cal E}_{d,r}$ with the 
desired properties. Since the quot scheme is a fine moduli space,
moreover, we know that $(\mbox{id},\beta _r)^*{\cal S}_d = {\cal E}_{d,r}$.

\medskip

It may be more illuminating to think of the maps $\beta _r$ pointwise.
Namely,
if $S \hookrightarrow V^*\otimes {\cal O}_{{\bf P}^1}$ is a vector bundle
subsheaf
of rank $k$ and degree $-d+r$, then a point $p\in {\bf P}^1$
and rank $r$ quotient $S(p) \rightarrow {\bf C}^r(p)$ 
determine a point $x\in {\cal G}_{d,r}$. The kernel
of the map $S \rightarrow {\bf C}^r(p)$ is a new vector bundle
$E$ of rank $r$ and degree $-d$ which becomes a subsheaf of
$V^*\otimes {\cal O}_{{\bf P}^1}$ via its inclusion as a subsheaf of $S$. 
The resulting subheaf $E\hookrightarrow V^*\otimes {\cal O}_{{\bf P}^1}$
is the image $\beta _r(x)$.

\medskip

Now, suppose that $V^*\otimes {\cal O}_{{\bf P}^1} \rightarrow T$ is
a quotient with $\chi({\bf P}^1,T(m)) = \chi$, and that the 
rank of $T$ at $p\in {\bf P}^1$ is at least $n-k+r$. Then let $i:E
\hookrightarrow V^*\otimes {\cal O}_{{\bf P}^1}$ be the kernel, and
consider the dual map $i^*$. The fact that $T$ has rank $n-k+r$ at $p$
implies that at $p$, the map on fibers: $i^*(p): V(p) \rightarrow E^*(p)$
has a cokernel of rank at least $r$. Thus, we may choose 
a quotient $E^* \rightarrow {\bf C}^r(p)$ such that 
$i^*$ factorizes through the kernel, $S^*$, which proves (i).
Moreover, if the sheaf $T$ has rank exactly $n-k+r$ at exactly
one point $p\in {\bf P}^1$, then the bundle $S^*$ is uniquely determined,
which proves (ii) on the level of sets.

\medskip

To prove (ii) completely, we observe that the map $\beta _r$ may be inverted
on the image of $\pi_r^{-1}({\bf P}^1\times {\cal M}_{d-r})$ by 
globalizing the previous paragraph. Namely, on this image, the cokernel $N$
of the map $u:V\otimes {\cal O} \rightarrow {\cal S}^*_d$
is torsion, supported on a section $Z$ of ${\bf P}^1\times
{\cal Q}_d$ over ${\cal Q}_d$, and of rank $r$ on its support. 
The projection of $Z$ to ${\bf P}^1$,
kernel ${\cal E}^*$ of the map ${\cal S}^*_d \rightarrow N$ (which is a bundle!),
and the cokernel of the map ${\cal S}_d|_Z \rightarrow {\cal E}|_Z$ 
will give us the inverse to $\beta _r$.

\medskip

Finally, reconsider the maps
$V\otimes {\cal O} \rightarrow \pi^*{\cal S}^*_{d-r} \rightarrow {\cal E}^*
_{d,r} = (\mbox{id},\beta _r)^*{\cal S}^*_d$ of vector bundles on ${\bf P}^1\times 
{\cal G}_{d,r}$. The latter map is an isomorphism off of 
the preimage of $\Delta \times {\cal Q}_{d-r}$, and when restricted
to the preimage of $\Delta \times {\cal Q}_{d-r}$, it  factors
through $\hat S^*_{d-r}$. Thus the degeneracy locus where 
$V_{k+i-a_i} \otimes {\cal O} \rightarrow
{\cal E}^*_{d,r}$ has kernel of rank $i$ is the union of 
the same degeneracy loci for $\pi^*{\cal S}_{d-r}$ generically and
for $\hat S^*_{d-r}$ on the preimage of $\Delta \times {\cal Q}_{d-r}$. 
Since the rank of ${\hat S}^*_{d-r}$ is $k-r$,
we get (iii) when we restrict the degeneracy loci
to $p\times {\cal G}_{d,r}$.
 
\bigskip

As our first application of the structure theorem, we 
will prove the second moving lemma.

\newpage 

{\bf Proof of Moving Lemma 2.2A:} Note that the Lemma is identical
to Lemma 2.2 in case $d = 0$, and anyway it is easy in case $d = 0$ because
${\cal Q}_0 = {\cal M}_0 = G$.
We prove the Lemma in general by induction on
the degree.

\medskip

Notice first of all that the codimension of the intersection cannot
be {\it larger} than $A$ because each $\overline W_{\vec a}(p)$
has codimension at most $|\vec a|$, 
by \cite{F1}, Theorem 14.3(b).
Since Lemma 2.2 already takes care of the restriction to ${\cal M}_d$, 
it therefore suffices to show that  
$$\overline W_{{\vec a}_1}(p_1) \cap ... \cap \overline W_{{\vec a}_N}(p_N)
\cap {\cal B}_d$$
has codimension greater than $A$ in ${\cal Q}_d$.

By the structure theorem, it suffices to show that:
$$\cap_{j=1}^N \left\{ \pi_r^{-1}({\bf P}^1\times \overline W_{{\vec a}_j}(p_j)) \cup 
\widehat W_{{\vec a}_j - \vec r}(p_j)\right\}$$ 
has codimension greater than $A - (\mbox{dim}({\cal Q}_d) - \mbox{dim}
({\cal G}_{d,r}))$ in each ${\cal G}_{d,r}$. (In fact, it suffices to show
this for $r=1$, but we will need the other cases below.)

Since the points are {\it distinct} and $\widehat W_{{\vec a}-{\vec r}}(p)$ is concentrated
in $\pi_r^{-1}(p\times {\cal Q}_{d-r})$, it follows
that the only nonempty intersections admit one or zero occurrances of an
$\widehat W_{{\vec a}_j-\vec r}(p_j)$. Moreover, since we are assuming the Lemma
for lower degree, we find that the intersection 
$\cap_{j=1}^N\pi_r^{-1}({\bf P}^1
\times \overline W_{{\vec a}_j}(p_j))$ has codimension exactly $A$
in ${\cal G}_{d,r}$, and since 
$\mbox{dim}
({\cal G}_{d,r}) < \mbox{dim}({\cal Q}_d)$ (either by the structure theorem
or a dimension count),
we only have to prove (rearranging indices!) that intersections of the form:
$$(\dagger) \ \ \pi_r^{-1}({\bf P}^1 \times \cap _{j=1}^{N-1}\overline W_{{\vec a}_j}(p_j))
\cap \widehat W_{{\vec a}_{N}- \vec r}(p_N)$$
are of codimension greater than $A - (\mbox{dim}({\cal Q}_d) - 
\mbox{dim}({\cal G}_{d,r}))$ in ${\cal G}_{d,r}$.

Now consider the (largest) open subscheme $U_{d,r}(p_N) \subset \pi_r^{-1}(p_N\times
{\cal Q}_{d-r})$ over which the 
restriction of the map $V\otimes {\cal O}\rightarrow \hat {\cal S}^*_{d-r}$
to $p_N\times {\cal G}_{d,r}$ is surjective. This restriction
determines a map (which we may call evaluation at $p_N$)
$ev_{p_N}:U_{d,r}(p_N) \rightarrow G(n-k+r,n)$. By the same argument as Lemma 2.2,
one concludes that for any $Z\subset U_{d,r}(p_N)$, the intersection
$Z\cap \widehat W_{\vec a_N - \vec r}(p_N)$ has codimension at least
$|{\vec a}_N| - r(n-k)$ (and greater if $(a_N)_{n-k} - r < 0$) in $Z$.
If we let $Z$ be the (codimension $A - |\vec a_N|$) intersection of 
the ${\bf P}^1 \times \overline W_{{\vec a}_j}(p_j)$ with $U_{d,r}(p_N)$,
then
the open subset of $(\dagger)$ obtained by restricting to 
$U_{d,r}(p_N)$ has codimension at least $A - r(n-k) + 1$ in
${\cal G}_{d,r}$, and from the dimension count:
$$\begin{array}{rcl}\mbox{dim}({\cal Q}_d) - \mbox{dim}({\cal G}_{d,r}) & =
& dn - [(d-r)n + 1 + r(k-r)] \\
& = & r(n-k) + r^2 - 1 \\ \end{array}$$
we get the desired result for the restriction of
$(\dagger)$ to $U_{d,r}(p_N)$.

On the other hand, by (i) of Theorem 3.1, if $x \in {\cal Q}_d$ is 
in the image of $\pi_r^{-1}(p_N\times {\cal Q}_{d-r})$ 
but not in the image of $U_{d,r}(p_N)$ and $r < k$,
then ${\cal T}_d$ has rank at least $n-k+r+1$ at $(p_N,x)$, so $x$ is in 
the image of ${\cal G}_{d,r+1}$. 

Recall that we needed to prove the codimension estimate for $(\dagger)$
on ${\cal G}_{d,1}$ (since this surjects birationally onto the boundary). 
We could get the estimate for the open intersection with $U_{d,1}(p_N)$,
and observed that the complement maps to the image of ${\cal G}_{d,2}$
(which is birational to ${\cal G}_{d,2})$.
By the same reasoning and induction
on $r$, we are therefore reduced to proving the codimension estimate 
for ${\cal G}_{d,k}$.
(In other words, we still need to consider
the case where $p_N$ is a base point.)

But in this case, we have: 
$$\begin{array}{lcl} {\cal G}_{d,k} & = & {\bf P}^1\times {\cal Q}_{d-k},\\
\widehat W_{{\vec a_N} - \vec k}(p_N) & = & p_N \times {\cal Q}_{d-k}\\
\end{array}$$
and the sum $\sum_{j=1}^{N-1}|\vec a_j| = A - |\vec a_N|$
is certainly at least $A - k(n-k)$. This implies that the codimension of
$(\dagger)$ in ${\cal G}_{d,k}$ is at least $A - k(n-k) + 1$, by induction on
the degree, 
and we obtain the last case by 
the same dimension count as before with $r = k$.

\bigskip

{\bf Remark 3.2:} If $\vec b = (b_1,...,b_{n-k+1})$ is an $(n-k+1)$-tuple of
integers satisfying $k \ge b_1 \ge ... \ge b_{n-k+1} \ge 0$, then
we define $\overline W_{\vec b}(o)$ as in 2.1A.
If $b_{n-k+1} = 0$, then $\overline W_{\vec b}(o) = 
\overline W_{\vec b^t}(o)$, where $\vec b^t = (b_1,...,b_{n-k})$.
However if $b_{n-k+1} \ne 0$, then $\overline W_{\vec b}(o) \subset {\cal B}_d$
and Theorem 3.1(iii) applies here, too, to give:
$$\beta_1^{-1}(\overline W_{\vec b}(o)) = 
\pi_1^{-1}({\bf P}^1 \times \overline W_{\vec b}(o))
\cup \widehat W_{\vec b - \vec 1}(o)$$

The proof of Moving Lemma 2.2A can be applied to 
distinct points $o,p_1,...,p_N \in {\bf P}^1$,
$\overline W_{\vec b}(o)$ ($b_{n-k+1}\ne 0$), and ``ordinary''
Schubert varieties
$\overline W_{\vec a_1}(p_1),...,\overline W_{\vec a_N}(p_N)$. 
In this case, It tells us that the 
intersection has the expected
codimension $|\vec b| + \sum|\vec a_i|$ in ${\cal Q}_d$,
and that the image under $\beta_1$ of:
$$\widehat W_{\vec b - \vec 1}(o) \cap
\pi_1^{-1}(\{o\}\times \overline W_{\vec a_1}(p_1) \cap ... \cap \overline 
W_{\vec a_N}(p_N) \cap {\cal M}_{d-1})$$
is Zariski dense in that intersection. But  
$\beta_1$ is an embedding 
when restricted to $\pi_1^{-1}(\{o\}\times {\cal M}_{d-1})$, 
so when the 
intersection consists of distinct points, they may be counted
in ${\cal Q}_d$ or in 
$\pi_1^{-1}(\{o\}\times {\cal M}_{d-1})$,
or even in $\pi_1^{-1}(\{o\}\times {\cal Q}_{d-1})$ (any extra points in
the intersection 
$\widehat W_{\vec b - \vec 1}(o) \cap
\pi_1^{-1}(\{o\}\times \overline W_{\vec a_1}(p_1) \cap ... \cap \overline 
W_{\vec a_N}(p_N))$ 
would map to extra intersection points in ${\cal Q}_d$).
 
\newpage

\nt {\bf 4. Quantum Schubert Calculus:} In this section, we use 
Theorem 3.1 to prove the quantum versions of Giambelli
and Pieri as stated in \S 1.

\bigskip

{\bf Proof of Quantum Giambelli:} Suppose $M(W_{\vec *}) = 
c\prod _{\vec a} W_{\vec a}^{n_{\vec a}}$ is some monomial 
in the Schubert varieties. We define the Gromov-Witten invariants
of $M$ in the obvious way, by setting 
$$\langle M(W_{\vec a}) \rangle_d := c\langle ...,W_{\vec a},...,W_{\vec a},...\rangle_d.$$ 
where each $W_{\vec a}$ appears $n_{\vec a}$ times on the right. This
definition extends in the obvious way to define Gromov-Witten invariants
of any collection $P_1(W_{\vec *}),...,P_N(W_{\vec *})$ of polynomials
in the Schubert varieties. Also, the Conclusion following Corollary 2.4
applies to show that the Gromov-Witten invariant defined in this way 
coincides with the degree of the product of the $P_i(\sigma_{\vec *})$,
thought of as polynomials in the generalized Schubert cohomology classes
(see Corollary 2.3)
when evaluated in the cohomology ring of the quot schemes ${\cal Q}_d$.

\medskip

Thus, for example, the Giambelli determinantal formula 
for $W_{\vec a}$, evaluated with a quantum product, becomes:
$$\sum_{d \ge 0} q^d(\sum_{\vec b} \langle W_{\vec b},\Delta_{\vec a}(W_*)
\rangle_d \ \sigma_{\vec b}^c).$$

If we put this together with Lemma 2.5, then the quantum Giambelli 
formula is equivalent to the statement:
$$\langle W_{\vec b},W_{\vec a}\rangle_d = \langle W_{\vec b}, \Delta_{\vec a}(W_*)\rangle_d
\ \mbox{for all $d \ge 0$ and Schubert classes $W_{\vec b}$}.$$

But because these invariants are just the evaluations of the corresponding
cohomology classes in the corresponding quot scheme, quantum Giambelli follows from the 
(even stronger!) assertion:  
$$\sigma _{\vec a} = \Delta_{\vec a}(\sigma_*)\ \mbox{in the cohomology ring of
each ${\cal Q}_d$}.$$

\medskip

Choose $p\in {\bf P}^1$. Then 
since $\sigma_{\vec a}$ is the image in cohomology of $\overline W_{\vec a}(p)$,
which by the moving lemma has  pure codimension $|\vec a|$ in ${\cal Q}_d$, this
statement is a direct application of a theorem of Kempf and 
Laksov  (\cite{KL}) to the universal map
$V\otimes {\cal O} \rightarrow  {\cal S}_d^*$ of vector bundles on 
${\bf P}^1 \times {\cal Q}_d$, or rather, to the restriction of the 
universal map
to $\{p\}\times {\cal Q}_d$.

\bigskip

{\bf Proof of Quantum Pieri:} There
is a {\it polynomial} identity:
$$\sigma_a  \Delta_{\vec a}(\sigma_*) = 
\sum_{\vec b} \Delta_{\vec b}(\sigma_*)$$
where $\Delta_{\vec a}$ and the $\Delta_{\vec b}$ are the Giambelli determinants
and $\vec b$ varies over all $(n-k+1)$-tuples $(b_1,...,b_{n-k+1})$
with $k \ge b_1 \ge a_1 \ge b_2 \ge ... \ge b_{n-k+1} \ge 0$.
(See Lemma A.9.4 of \cite{F1} for a proof of this.)

Note that the $\vec b$ are not $(n-k)$-tuples! In the classical Schubert 
calculus, any $\vec b$ with nonzero $b_{n-k+1}$ gives rise to 
the empty variety in the Grassmannian, so if one sets those 
$\sigma_{\vec b}$ to zero,
then the classical Pieri's formula results. However, when we 
evaluate them in the cohomology ring of ${\cal Q}_d$ for 
positive $d$, we have seen in Remark 3.2 that
such $\vec b$ may give rise to nonzero varieties. In fact, I claim that
if $d > 0$ and $W_{\vec a}$ is any Schubert variety, then:
$$\langle \Delta_{\vec b}(W_*), W_{\vec a}\rangle _d = 
\left\{ \begin{array}{l} 0 \ \mbox{if $b_{n-k+1} > 0$ and $b_1 < k$} \\
\langle W_{(b_2-1,...,b_{n-k+1}-1)},W_{\vec a}\rangle_{d-1} \ \mbox{if 
$b_{n-k+1} > 0$ and $b_1 = k$} \end{array}\right.$$

We'll prove this claim later. Let's first see how quantum Pieri follows. 

\medskip

The polynomial identity above, together with
quantum Giambelli and classical Pieri, 
gives the following identity among quantum products:
$$(\dag)\ \sigma_a * \sigma_{\vec a} =
\sigma_a * \Delta_{\vec a}(\sigma_*) = 
p_{a,\vec a}(\sigma_{\vec *}) + 
\sum _{\{\vec b\ |\ b_{n-k+1} \ne 0\}} \Delta_{\vec b}(\sigma_*)$$

Multiply the formula in Lemma 2.5 by $q$ to get:
$$q\sigma_{\vec c} =  \sum_{d > 0} q^d (\sum_{\vec a}  \langle W_{\vec a}, 
W_{\vec c}\rangle _{d-1} \ \sigma_{\vec a^c})$$
for any Schubert cohomology class $\sigma_{\vec c}$. Using this formula applied
to the Schubert cohomology class $\sigma_{(b_2 - 1,...,b_{n-k+1}-1)}$ together
with the claim,
we see that the last terms in $(\dag)$ evaluate under the quantum product as
follows:
$$\Delta_{\vec b}(\sigma_*) = \left\{
\begin{array}{l} 0 \ \mbox{if $b_1 < k$ and} \\
q\sigma_{(b_2-1,...,b_{n-k+1}-1)} \ \mbox{if}\ b_1 = k
\end{array}\right.$$

Putting this together with $(\dag)$ gives the quantum Pieri formula.

\medskip

{\bf Proof of the claim:} Suppose that $\vec b = (b_1,...,b_{n-k+1})$
satisfies the inequalities $k \ge b_1 \ge b_2 ... \ge b_{n-k+1} > 0$. 
Then by the theorem of Kempf-Laksov again, together with Remark 3.2,
we know that when evaluated in the cohomology ring of ${\cal Q}_d$,
the Giambelli determinant $\Delta_{\vec b}(\sigma_*)$ is equal to 
$\sigma_{\vec b}$, the image in cohomology 
of the degeneracy locus $\overline W_{\vec b}(p)$.

\medskip

For such $\vec b$, let's {\bf define} the Gromov-Witten invariants
$\langle W_{\vec b},W_{\vec a}\rangle_d$ to be the degree in the 
cohomology ring of ${\cal Q}_d$ of the product of $\sigma_{\vec b}$ and
$\sigma _{\vec a}$. Equivalently, these Gromov-Witten invariants are
the number of points 
{\bf in the quot scheme} (as opposed to ${\cal M}_d$, which would trivially
give zero) in 
$\overline W_{\vec b}(o) \cap \overline W_{\vec a}(p)$ for general translates
of the flags and distinct points $o,p \in {\bf P}^1$. Then with this definition, 
the claim is equivalent to the equalities:
$$\langle W_{\vec b}, W_{\vec a}\rangle _d = 
\left\{ \begin{array}{l} 0 \ \mbox{if $b_{n-k+1} > 0$ and $b_1 < k$} \\
\langle W_{(b_2-1,...,b_{n-k+1}-1)},W_{\vec a}\rangle_{d-1} \ \mbox{if 
$b_{n-k+1} > 0$ and $b_1 = k$} \end{array}\right.$$
for all $d > 0$ and $W_{\vec a}$.

\medskip

Notice that the dimensions work out(!) In other words, if $b_1 = k$, then: 
$$|\vec b| + |\vec a| = \ \mbox{dim}({\cal Q}_d) 
\Leftrightarrow |(b_2-1,...,b_{n-k+1}-1)| + |\vec a| = \ \mbox{dim}({\cal Q}_{d-1})$$

We assume that this equality holds (otherwise the claim is trivial). Then 
by Remark 3.2,
$$\langle W_{\vec b},W_{\vec a}\rangle_d =  
\ \mbox{number of points in}\ \ 
\widehat W_{\vec b - \vec 1}(o) \cap \pi_1^{-1}(o\times \overline W_{\vec a}(p)) $$
where $\pi_1^{-1}(\{o\} \times {\cal Q}_{d-1}) \subset {\cal G}_{d,1}$
is the projectivization of the restriction ${\cal S}_d(o)$ 
of ${\cal S}_d$
to $\{o\}\times {\cal Q}_{d-1}$, and $\widehat W_{\vec b - \vec 1}(o)$ 
is the degeneracy locus for the map from $V\otimes {\cal O}$ to 
$\widehat {\cal S}^*$, where $\widehat {\cal S}$ is the universal subbundle
of $\pi_1^*{\cal S}_d(o)$. 

\medskip

Thus the claim follows if we can show that in the homology of 
the smooth varieties 
$\pi_1^{-1}(\{o\}\times {\cal Q}_{d-1}) \subset {\cal G}_{d,1}$ and ${\cal Q}_{d-1}$, 
$$(\pi_1)_*([\widehat W_{\vec b - \vec 1}(o)].\pi_1^*[\overline W_{\vec a}(p)]) = 
\left \{ \begin{array}{l} 0 \ \mbox{if} \ b_1 < k \ \mbox{and} \\
\left[\overline W_{(b_2-1,...,b_{n-k+1}-1)}(o)\right].[\overline W_{\vec a}(p)] \ \mbox{if}\ b_1 = k\end{array} \right.$$

But Kempf-Laksov applied to $\widehat W_{\vec b - \vec 1}(o)$ and
$\overline W_{(b_2 - 1,...)}(o)$ makes this a special case 
of a formula of J\'osefiak, Lascoux and Pragacz (see Example 14.2.2 
of \cite{F1} and \cite{JLP}).

\newpage

\nt {\bf 5. The Formula of Vafa and Intriligator:} There is a marvelous
residue formula due to Vafa and Intriligator which uses the presentation 
$(*)_q$ of the quantum cohomology ring to compute the
Gromov-Witten intersection numbers:
$$\langle W_{a_1},...,W_{a_N}\rangle _d$$
of {\bf special} Schubert varieties. (See \cite{I},\cite{ST}, \cite{B}.)
The formula is the following:

\bigskip

\nt {\bf (Vafa and Intriligator's)  Formula:} Fix  
$\zeta$ a primitive $n$th root of $(-1)^k$ and assume that $0 \le a_i \le k$ and
$a_1+ ... + a_N = \ \mbox{dim}({\cal M}_d)$.  Then:

\medskip

\nt $\langle W_{a_1},...,W_{a_N}\rangle_d = $
$$(-1)^{\left( k \atop 2 \right)}n^{-k}\sum_{i_1 > ... > i_k} \sigma_{a_1}(\zeta^I)
\cdots \sigma_{a_N}(\zeta^I)
\left(\frac {\prod_{j \ne l}(\zeta^{i_j}-\zeta^{i_l})}
{\prod_{j=1}^k\zeta^{(n-1)i_j}}\right)$$
where $\zeta^I = (\zeta^{i_1},...,\zeta^{i_k})$ and $\sigma_{a_i}$ are the elemetary
symmetric polynomials in $k$ variables (i.e. $\sigma_0(\zeta^I) = 1, 
\sigma_1(\zeta^I) = \zeta^{i_1} + ... + \zeta^{i_k}$, etc.)

\medskip

The point I want to make is that because of quantum Giambelli, 
the same formula computes all the Gromov-Witten intersection numbers .
That is:

\medskip

\nt {\bf Corollary (of quantum Giambelli):} Assume $W_{\vec a_1},...,W_{\vec a_N}$
are Schubert varieties on $G$ satisfying $|\vec a_1| + ... + |\vec a_N| = 
\ \mbox{dim}({\cal M}_d)$. Then the Gromov-Witten intersection 
number:
$$\langle W_{\vec a_1},...,W_{\vec a_N}\rangle _d$$
may be computed by the Vafa-Intriligator formula, where
the elementary symmetric polynomials $\sigma_{a_i}(\zeta^I)$
are replaced by the Giambelli determinants $\Delta_{\vec a_i}(\sigma_*(\zeta^I))$ of the 
elementary symmetric polynomials.

\bigskip

University of Utah, Salt Lake City, UT 84112

\medskip

email: bertram@math.utah.edu

\end{document}